# BUILDING AND USING SEMIPARAMETRIC TOLERANCE REGIONS FOR PARAMETRIC MULTINOMIAL MODELS

By Jiawei Liu and Bruce G. Lindsay[1]

*Georgia State University and Pennsylvania State University*

We introduce a semiparametric "tubular neighborhood" of a parametric model in the multinomial setting. It consists of all multinomial distributions lying in a distance-based neighborhood of the parametric model of interest. Fitting such a tubular model allows one to use a parametric model while treating it as an approximation to the true distribution. In this paper, the Kullback–Leibler distance is used to build the tubular region. Based on this idea one can define the distance between the true multinomial distribution and the parametric model to be the index of fit. The paper develops a likelihood ratio test procedure for testing the magnitude of the index. A semiparametric bootstrap method is implemented to better approximate the distribution of the LRT statistic. The approximation permits more accurate construction of a lower confidence limit for the model fitting index.

**1. Introduction.** The conventional approach for assessing goodness of fit in multinomial models is discussed in several standard sources [see, e.g., Bishop, Fienberg and Holland (1975) and Agresti (1990)]. The agreement between the model and the data is assessed with a goodness-of-fit test statistic, such as the Pearson chi-squared statistic, the likelihood ratio statistic or other measures of divergence [Read and Cressie (1988)], and the model is then accepted as true or rejected as false. These methods for evaluating rely on $\chi^2$-statistics or quantities derived from them. In cases where the sample size $n$ is sufficiently large, the model will usually be rejected in favor of some more complicated model. Goodness-of-fit statistics are usually not informative when the sample size is very large. The question they fail to address is: "Is the specified model a good approximation to the true distribution of the sample?"

Received May 2005; revised March 2008.
[1]Supported in part by NSF Grant DMS-04-05637.

In 1994, Rudas, Clogg and Lindsay presented a framework based on mixture methods for evaluating goodness of fit of contingency tables. Suppose $\boldsymbol{\tau}$ is the true cell probability. For a given parametric model $\mathcal{M}$, $\tau$ can generally be written as a two-point mixture as $\boldsymbol{\tau} = (1-\pi)\mathbf{m} + \pi\mathbf{e}$, where $\mathbf{m}$ is an element from the model $\mathcal{M}$, $\mathbf{e}$ is an unspecified distribution corresponding to the modeling error, and $\pi$ is the mixture weight. The *mixture index of fit*, $\pi^*$, is defined to be the smallest such $\pi$. That is, $\pi^*$ is the fraction of the population that could not possibly be described by model $\mathcal{M}$. The larger the value of $\pi^*$, the more discrepant the true distribution is from any model distribution. Thus one can think of $\pi^* = \pi^*(\tau, \mathcal{M})$ as measuring a distance from $\tau$ to the model $\mathcal{M}$. This approach focuses on measuring the difference between truth and model, and so on the substantive importance of the discrepancy between the model and the data. This permits an evaluation of the model that downplays the role of sample size.

In this paper, we extend the results of Rudas et al. to the Kullback–Leibler distance, $K^2(\boldsymbol{\tau}, \mathbf{m}) = \sum \mathbf{m} \log[\mathbf{m}/\boldsymbol{\tau}]$. Our motivation for choosing this distance is two fold. First, the $\pi^*$ distance measure has some undesirable features. We can write

$$\pi^*(\boldsymbol{\tau}, \mathcal{M}) = \min_\theta \max_t \left\{ 1 - \frac{\tau(t)}{m(t)} \right\},$$

where $t$ is the cell index and $\theta$ is the parameter of the model [Xi (1996), page 14]. The function of $\pi^*$ is not everywhere differentiable, which generates a non-standard asymptotic theory [Xi (1996), page 78] and difficulties in computation [Xi and Lindsay (1996)]. On the other hand $K^2$ is smooth and generates asymptotically efficient parameter estimators, as we will show. A reason for choosing $K^2$ from among smooth distances is the simple and useful geometric structure that arises from $K^2$, to be described below.

Given any statistical distance function $\rho(\boldsymbol{\tau}, \mathbf{m})$, such as $K^2(\boldsymbol{\tau}, \mathbf{m})$, we can define the statistical distance from the model $\mathcal{M}$ to the true density $\boldsymbol{\tau}$ via

(1) $$\rho(\boldsymbol{\tau}, \mathcal{M}) = \inf_{\mathbf{m} \in \mathcal{M}} \rho(\boldsymbol{\tau}, \mathbf{m}).$$

The index $\rho^* = \rho(\boldsymbol{\tau}, \mathcal{M})$ is then a measure of the distance from the parametric model $\mathcal{M}$ to the true distribution $\boldsymbol{\tau}$. We can then build a tolerance region about $\mathcal{M}$, which consists of all distributions $\boldsymbol{\tau}$ that satisfy $\rho(\boldsymbol{\tau}, \mathbf{m}) \leq c$, where $c$ is some specified level of tolerance.

In this paper we consider the problem of testing the hypothesis

(2) $$H_0 : K^2(\boldsymbol{\tau}, \mathcal{M}) \leq c.$$

That is, does the model $\mathcal{M}$ provide an acceptable approximation to true distribution $\boldsymbol{\tau}$, where the word "acceptable" is measured by $K^2$ and indexed by constant $c$. We will test the hypothesis using the likelihood ratio test statistic in Section 3.1, in common with Rudas, Clogg and Lindsay (1994). Our reason for using likelihood ratio test (LRT) is the well-known power of this procedure in detecting deviations from the null hypothesis. Of course, given a test statistic for each $c$-level, one can invert the LRT to obtain a lower confidence limit for $K^2(\boldsymbol{\tau}, \mathcal{M})$.

An interesting and very useful feature arises from using the likelihood ratio as a test statistic together with $K^2(\boldsymbol{\tau}, \mathbf{m})$ as the distance. We show in Section 3 that the maximum likelihood estimator $\hat{\mathbf{p}}$ of the probability distribution under $H_0$ has the form

$$(3) \qquad \hat{\mathbf{p}} = \pi \mathbf{d} + (1-\pi) \mathbf{m}_{\hat{\theta}},$$

where $\mathbf{m}_{\hat{\theta}}$ is an estimator of the model element closest to $\boldsymbol{\tau}$, $\mathbf{d}$ is the vector of observed data and the weight $\pi \in [0,1]$ depends on $c$. That is, viewed geometrically, the null estimator is a convex combination mixture of the best model distribution $\mathbf{m}_{\hat{\theta}}$ and the data proportions $\mathbf{d}$. This leads to a very simple method, essentially an iteratively reweighted maximum likelihood algorithm, for finding the maximum likelihood estimator $\hat{\theta}$ and hence the LRT statistic. Section 3 concludes by showing why the estimator $\hat{\theta}$ is both efficient under the model and robust under model failure.

In Section 4 we consider two methods to set critical values for the test. One is based on asymptotic distribution theory and the other, more accurate, method is based on a semiparametric bootstrap using (3), the null hypothesis estimator. Our first real data example is in Section 5, where we compare our methods with the conventional likelihood method on two simple data sets. Section 6 contains the analysis of a richer multiway table, and we conclude in Section 7.

**2. Distance-based tubular models.** Given $\rho(\boldsymbol{\tau}, \mathbf{m})$, a distance, and $\rho(\boldsymbol{\tau}, \mathcal{M})$ defined in (1), we will say that the model $\mathcal{M}$ is *adequate at level* $c$ if $\rho(\boldsymbol{\tau}, \mathcal{M}) \leq c$. This is equivalent to saying that $\boldsymbol{\tau}$ is in the model tube $\mathcal{M}_c$, where the *model tube* $\mathcal{M}_c$ is defined as all multinomial densities $\mathbf{p}$ sufficiently close to the model $\mathcal{M}$:

$$\mathcal{M}_c = \bigcup_{\mathbf{m} \in \mathcal{M}} \{\mathbf{p} : \rho(\mathbf{p}, \mathbf{m}) \leq c\}.$$

This representation shows that a tube can be described as the union of balls around each model element. Since $\boldsymbol{\tau} \in \mathcal{M}_c \Leftrightarrow \rho(\boldsymbol{\tau}, \mathcal{M}) \leq c$, we can think of $\mathcal{M}_c$ as the extended model defined by the null hypothesis (2).

2.1. *Selecting the tolerance level* $c$. The idea of using a tolerance zone around a parametric model for goodness-of-fit testing can be found in Hodges and Lehmann (1954), Goutis and Robert (1998) and Dette and Munk (2003), as well as Rudas, Clogg and Lindsay (1994).

One important and challenging issue with the tubular hypothesis is the choice of the bound $c$. Rudas, Clogg and Lindsay (1994) supposed the true distribution can be written as a two-point mixture $\boldsymbol{\tau} = (1-\pi)\mathbf{m} + \pi\mathbf{e}$. The $\pi^*$ index is defined as the smallest such $\pi$. The $\pi^*$ index has a natural mixture interpretation, so $c = 0.05$, say, can be thought as 5% contamination. However, if one uses a standard distance in the other three papers, it is very difficult to give a direct interpretation to $c = 0.05$.

Hodges and Lehmann (1954) proposed the tolerance zone around the null hypothesis in a number of testing problems. They used ordinary Euclidean distance or a weighted Euclidean distance as the distance from a model element to the truth, an example, the Neyman chi-squared distance. Hodges and Lehmann did not give a detailed discussion on how one should choose $c$.

Goutis and Robert (1998) proposed a Bayesian model selection approach using a tubular model based on the likelihood distance. The problem to generate a proper bound $c$ had been considered by Mengersen and Robert (1996) and Dupuis (1997). In their models the distance was bounded by a constant. They selected $c$ to be a fixed fraction of this bound. However, the choice of the fixed fraction is still an issue.

Dette and Munk (2003) discussed the tubular hypothesis with ordinary Euclidean distance for nonparametric regression models. They thought that the choice of $c$ was one of the "most difficult tasks." Dette and Munk considered model selection more as an explorative data analysis and let $c = 0$, 0.5 and 1 in their examples. Dette and Munk suggested that a proper analysis of the power would be helpful.

Although there is no perfect way to interpret a standard distance in the literature, we do offer two additional assessment tools that can provide guidance in deciding on an acceptable distance. We used the Kullback–Leibler distance to build the tubular model.

First, the square root of $K^2(\boldsymbol{\tau}, \mathcal{M})$, is similar in magnitude to the $\pi^*$ index of Rudas, Clogg and Lindsay (1994), as we found in examples of Sections 5 and 6. This similarity helps us to interpret the distance. The index $\pi^*$ is the proportion of the population outside of the model. So the value of root $K^2$ roughly measures the same degree of faith on the model, with $\pi^* = 0$ ($\sqrt{K^2}$ close to zero) representing full use of the model estimation and $\pi^* = 1$ ($\sqrt{K^2}$ close to 1) representing complete discard of the model.

Secondly, another very statistical way to think about the distance between two distributions $F$ and $G$ is to ask about the sampling consequences of the distance. In particular, how different are samples of size $n$ from these two distributions? One could, for example, ask what the power would be if one were to take one sample from each and test the null hypothesis that they were from the same distribution. This idea is due to Davies (1995). The problem with this measure is that it is dependent on the sample size $n$, and so not an absolute measure of difference.

Lindsay and Liu (2005) proposed an interpretable measure by determining the sample size $N^*$ at which one would have 50% power in a test of $F$ versus $G$, with the idea that at 50% power, the differences in the samples between $F$ and $G$ are starting to become obvious. Choosing the power 50% is convenient for calculation. It enables us to avoid calculating the variance of the asymptotic distribution [Lindsay and Liu (2005)]. We will use this secondary measure $N^*$ as an aid in the interpretation of radius $c$ in our data sections.

**3. Likelihood ratio test for tubular model.** In this section we will construct the likelihood ratio test for the hypotheses $H_0 : \boldsymbol{\tau} \in \mathcal{M}_c$ vs. $H_1 : \boldsymbol{\tau} \notin \mathcal{M}_c$. The test statistic is described, and an iteratively reweighted algorithm is developed to find the null estimators. Interpretations of the estimators are discussed at the end of this section.

3.1. *Finding the test statistic.* We will use the likelihood ratio test (LRT) for testing the null hypothesis (2). When the model is multinomial with the observed proportions $\mathbf{d}$, the *best model element* in tubular model $\mathcal{M}_c$ is defined by

$$\hat{\mathbf{p}}_c = \arg \min_{\mathbf{p} \in \mathcal{M}_c} L^2(\mathbf{d}, \mathbf{p}),$$

where $L^2$ is the likelihood distance, namely $L^2(\mathbf{d}, \mathbf{p}) = \sum d(t) \log(d(t)/p(t))$. The best model element is also the maximum likelihood estimator of the true multinomial probability density $\boldsymbol{\tau}$ under the null. The statistic $2nL^2(\mathbf{d}, \mathcal{M}_c) = 2nL^2(\mathbf{d}, \hat{\mathbf{p}}_c)$ equals the LRT statistic for testing the hypothesis $H_0 : \boldsymbol{\tau} \in \mathcal{M}_c$ vs. $H_1 : \boldsymbol{\tau} \notin \mathcal{M}_c$.

The problem of finding the test statistic $L^2(\mathbf{d}, \hat{\mathbf{p}}_c)$ can be written as the following optimization problem:

(4) $$\min_{\mathbf{p}} L^2(\mathbf{d}, \mathbf{p}) \quad \text{subject to} \quad \mathbf{p} \in \mathcal{M}_c.$$

That is, we minimize the likelihood distance subject to the constraint $K^2(\mathbf{p}, \mathcal{M}) \leq c$.

In this section, we show how the special structure of $L^2$ and $K^2$ can be used to turn (4) into a simple unconstrained problem. To do so, we first consider a simpler problem. Suppose that $\mathcal{M}$ contains a single density, that is, $\mathcal{M} = \{\mathbf{m}_0\}$. The simpler optimization problem is then

(5) $$\min_{\mathbf{p}} L^2(\mathbf{d}, \mathbf{p}) \quad \text{subject to} \quad K^2(\mathbf{p}, \mathbf{m}_0) \leq c.$$

THEOREM 1. *The solution $\hat{\mathbf{p}}_c$ for the optimization problem (5) exists on the boundary of $\mathcal{M}_c$ and has the form*

$$(6) \qquad \hat{\mathbf{p}}_c = \pi \mathbf{d} + (1-\pi)\mathbf{m}_0,$$

*where $\pi$ is a unique value in $[0,1]$ forcing the solution $\hat{\mathbf{p}}_c$ to be on the tube boundary, so $K^2(\hat{\mathbf{p}}_c, \mathbf{m}_0) = c$.*

PROOF. We use the method of Lagrange multipliers for constraints $K^2(\mathbf{p}, \mathbf{m}_0) \leq c$ and $\sum p(t) = 1$.  □

That is, $\hat{\mathbf{p}}_c$ is on the line connecting $\mathbf{m}_0$ and $\mathbf{d}$ and is on the boundary of the tube. It is interesting to note that if we replaced both $L^2$ and $K^2$ with the Euclidean distance, $\sum [d(t) - m(t)]^2$, we would have the same answer as (6). Heuristically, the spatial curvatures of $K^2$ and $L^2$ seem to cancel each other out, yielding a Euclidean-like structure.

Note that if we choose a fixed $\pi$-value and solve the unconstrained problem

$$(7) \qquad \min_{\mathbf{p}} \pi L^2(\mathbf{d}, \mathbf{p}) + (1-\pi)K^2(\mathbf{p}, \mathbf{m}_0)$$

then the solution $\hat{\mathbf{p}}_\pi$ has the same form as (6). Observe that the constraints on $c$ and $\pi$ are different in problems (5) and (7). In (5), $c$ is a preset constant and $\pi$ is a free parameter that must be chosen to drag the solution to the tube boundary and hence is a function of $c$. But in (7), $\pi$ is a pre-selected constant but it corresponds to a value of $c$, via

$$(8) \qquad c = K^2(\hat{\mathbf{p}}_\pi, \mathbf{m}_0),$$

so $c$ is a function of $\pi$. Hence there is a one-to-one relationship between $\pi$ and $c$. This suggests a strategy for finding the tubal test statistic for a fixed $c$. The optimization problem (7) is easier to solve. If we solve (7) on a dense grid of $\pi$ values, they will generate a set of corresponding $c$ values using (8), some of which would be quite close to the radius we want.

3.2. *Optimization steps and an algorithm.* Suppose $\mathcal{M}$ is not a single element $\{\mathbf{m}_0\}$. Problem (7) suggests consideration of the extended problem

$$(9) \qquad \min_{\mathbf{p}} \pi L^2(\mathbf{d}, \mathbf{p}) + (1-\pi)K^2(\mathbf{p}, \mathcal{M}).$$

Remarkably, once again a solution $\hat{\boldsymbol{p}}_\pi$ to problem (9) for fixed $\pi$ is also a solution $\hat{\mathbf{p}}_c$ to problem (4) for the appropriate choice of $c$, namely $K^2(\hat{\mathbf{p}}_\pi, \mathcal{M})$ [Liu (2003)].

For models $\mathcal{M} = \{\mathbf{m}_\theta\}$ with more than one element, (9) is equivalent to

$$(10) \qquad \min_{\theta} \min_{\mathbf{p}} \pi L^2(\mathbf{d}, \mathbf{p}) + (1-\pi)K^2(\mathbf{p}, \mathbf{m}_\theta).$$

Then for fixed $\theta$, the inner minimization in (10) has the simple explicit solution $\hat{\mathbf{p}}_\theta = \pi\mathbf{d} + (1-\pi)\mathbf{m}_\theta$. This ensures that the final solution is a simple convex combination of $\mathbf{d}$ and a model element $\mathbf{m}_{\hat\theta}$.

Next, to minimize over $\theta$, the optimization problem becomes

$$\min_\theta \pi L^2(\mathbf{d}, \hat{\mathbf{p}}_\theta) + (1-\pi)K^2(\hat{\mathbf{p}}_\theta, \mathbf{m}_\theta). \tag{11}$$

This minimization problem is one in a class of minimum distance problems considered in Lindsay (1994). Basu and Lindsay (2003) introduced an iteratively reweighted estimation function approach. We apply the results.

The estimating equation of (11) can be written as a weighted form of likelihood function

$$\sum \omega(\delta(t))d(t)u_\theta(t) = 0, \tag{12}$$

where $\delta$ is the *Pearson Residual Function* $\delta(t) = (d(t) - m_\theta(t))/m_\theta$, $\mathbf{u}_\theta = \nabla_\theta \mathbf{m}_\theta / \mathbf{m}_\theta$, and the weights $\omega = \log[1+\pi(\delta+1)/(1-\pi)]/(1+\delta)$.

The algorithm alternates the following steps until convergence.

- Given current estimate $\theta$, create weight $\omega = \log[1+\pi(\delta(t)+1)/(1-\pi)]/(1+\delta(t))$.
- Solve for the new estimate of $\theta$ from the estimating equation (12).

Note when $c = 0$, the maximum likelihood estimating equation for the original model is

$$\sum d(t)u_\theta(t) = 0. \tag{13}$$

Compare (13) with the estimating equation (12) of the tubular model. We can think of solving (12) as obtaining the maximum likelihood estimate from "pseudo data" $d^*(t) = \omega(\delta(t))d(t)$. The advantage is that for fixed $\omega(\delta(t))$ one can easily implement standard maximum likelihood statistical software to solve this reweighted estimating equation.

3.3. *Interpretation of the estimates.* For a second interpretation of $\hat{\mathbf{p}}_\pi$, note that an equivalent objective function to (9) is $L^2(\mathbf{d}, \mathbf{p}) + \gamma K^2(\mathbf{p}, \mathcal{M})$, for $\gamma = (1-\pi)/\pi$. It has the same solution $\hat{\mathbf{p}}_\pi$ as the maximum of the penalized likelihood function

$$P_\gamma(\mathbf{p}) = \sum d(t)\ln p(t) - \gamma K^2(\mathbf{p}, \mathcal{M}),$$

where $\gamma > 0$ is the penalty parameter. This criterion gives an estimator $\hat{\mathbf{p}}_\gamma$ that makes the likelihood large while being penalized, for being distant from the model $\mathcal{M}$ in $K^2$ terms.

Next, we consider interpretation of the estimators $\hat\theta = \hat\theta_\pi$ obtained from (11). If we define a new distance

$$T_\pi^2 = \pi L^2(\mathbf{d}, \pi\mathbf{d} + (1-\pi)\mathbf{m}_\theta) + (1-\pi)K^2(\pi\mathbf{d} + (1-\pi)\mathbf{m}_\theta, \mathbf{m}_\theta),$$

then $\hat{\theta}_\pi$ is the corresponding minimum distance estimator. When $\pi = 0.5$, this distance is symmetric in **d** and **m** and of some interest as a symmetric version of Kullback–Leibler distance. We call the symmetric distance the *mid-tube distance*, written as

$$T_{1/2}^2(\mathbf{d},\mathbf{m}) = \tfrac{1}{2}L^2(\mathbf{d},\mathbf{d}/2+\mathbf{m}/2) + \tfrac{1}{2}K^2(\mathbf{d}/2+\mathbf{m}/2,\mathbf{m}).$$

One can show, using the methods of Lindsay (1994), that the resulting estimators are, for all $\pi$, asymptotically as efficient as maximum likelihood when the model is correct. Moreover, they increase in robustness as $\pi$ decreases from 1 to 0, where $\pi = 1$ is maximum likelihood. Indeed, the mid-tube distance is topologically equivalent to Hellinger distance, well known for its robustness [Liu (2003)].

The mid-tube distance has the advantage in case that there is a cell with zero $\boldsymbol{\tau}$-probability and positive **m**-probability (or zero **m** and positive $\boldsymbol{\tau}$), where some conventional statistics, for example Pearson chi-square, would fail.

**4. Distribution theory.** To determine whether to reject the null hypothesis, we need to know the limiting distribution of the test statistic. In this section, we approximate the distribution of the test statistics $2nL^2$ by two methods. One is based on asymptotic distribution theory and the other method is based on a semiparametric bootstrap using the null hypothesis estimator (3). The latter one is found to be more accurate, especially for small values of $n$ and $c$.

4.1. *Asymptotic distribution of LRT statistics.* The asymptotic distribution of our likelihood ratio statistic $2nL^2$ depends on $c$. If $c = 0$, we are exactly in the setting of the standard likelihood ratio test of model $\mathcal{M}$ against the unrestricted multinomial alternative. And so under model regularity, the likelihood ratio test statistic $2nL^2$ has an asymptotic $\chi^2$ distribution with degrees of freedom equal to the number of cells $N$ minus the number of nonredundant parameters in $\mathcal{M}$, minus 1. This is the standard type of limiting distribution for likelihood ratio testing.

However, if $c > 0$, the LRT statistic $2nL^2$ has, asymptotically, a mixed $\chi^2$ distribution, with probability 0.5 equal to $\chi_0^2$ (point mass at zero) and with probability 0.5 equal to $\chi_1^2$. That is

(14) $$2nL^2 \xrightarrow{\mathcal{D}} \tfrac{1}{2}\chi_0^2 + \tfrac{1}{2}\chi_1^2.$$

This result was given for chi-squared distance in Hodges and Lehmann (1954) and for $\pi^*$ in Rudas, Clogg and Lindsay (1994).

If we use this asymptotic result (14), the test should reject the null hypothesis whenever the test statistic exceeds $\chi_1^2(0.90) = 2.70$, giving asymptotic size $\alpha = 0.05$. Inverting this test procedure, we can obtain a lower 95% confidence limit for the index $\rho^* = \rho(\mathbf{d}, \mathcal{M})$ defined in Section 1. It equals $\hat{\rho}_L^*$, where the tube radius $\hat{\rho}_L^*$ satisfies $2nL^2(\mathbf{d}, \mathcal{M}_{\hat{\rho}_L^*}) = 2.70$.

4.2. *Simulated distribution through bootstrapping.* One should be careful of using the asymptotic distribution in the cases that $c$ is almost zero and the sample size $n$ is not very large. In this section we offer an alternative method to simulate the distribution of $L^2$ using a semi-parametric bootstrap.

The logic of the mixture distribution is that for $c > 0$, the tube $\mathcal{M}_c$ has an open interior. When $\mathbf{d}$ falls inside the tube, the LRT statistic is zero. For a tube with a smooth boundary, asymptotically its tangent hyperplane can be a good approximation to the tube surface. If the data falls into the model side of the tangent hyperplane, the likelihood deviation is zero; and it is one-dimensional normal otherwise, so the squared deviation is $\chi_1^2$. The mixture weight is the probability that the data are on one side of a hyperplane, which is clearly 0.5.

The discontinuity in the limiting distribution at $c = 0$ arises because the dimension of $\mathcal{M}_0$ does not match that of $\mathcal{M}_c$ for $c > 0$. If $c$ is close to zero, $\mathcal{M}_c$ has almost no interior and so the LRT statistic is hardly zero. The limiting distribution would be better described by the chi-squared distribution as $c = 0$.

Moreover, in the multivariate case, when the sample size $n$ is not very large, the difference between the tangent hyperplane and the tube surface caused by the curvature of the tube can be significant. The mixture weight 0.5 is not reliable.

In this situation, it would be desirable to construct critical values from a method that takes into account the finite sample failure of the asymptotic approximations. We propose to use a semiparametric bootstrap method, which does just that.

Examples illustrating the implementation of these testing methods are in the next section.

TABLE 1
*Cross-classification of eye color and hair color ($n = 592$)*

| Eye color | Hair color | | | |
|---|---|---|---|---|
| | **Black** | **Brunette** | **Red** | **Blonde** |
| Brown | 68 | 119 | 26 | 7 |
| Blue | 20 | 84 | 17 | 94 |
| Hazel | 15 | 54 | 14 | 10 |
| Green | 5 | 29 | 14 | 16 |

TABLE 2
*Cross-classification of number of children by annual income ($n = 25,263$)*

| No. of children | Annual income | | | |
|---|---|---|---|---|
| | 0–1 | 1–2 | 2–3 | 3+ |
| 0 | 2161 | 3577 | 2184 | 1636 |
| 1 | 2755 | 5081 | 2222 | 1052 |
| 2 | 936 | 1753 | 640 | 306 |
| 3 | 225 | 419 | 96 | 38 |
| 4+ | 39 | 98 | 31 | 14 |

**5. Application to the independence model.** In this section we give a detailed examination of the tubular model in two-way contingency tables. The model $\mathcal{M}$ will be the row–column independence model.

5.1. *Two well-known examples.* The following examples were used by Diaconis and Efron (1985) to show problems associated with the conventional statistics for evaluating the model adequacy.

For Table 1 the LRT statistic is 146.44 on 9 degrees of freedom. The model would be rejected on the basis of these quantities. In 1994, Rudas, Clogg and Lindsay calculated the mixture index of fit $\pi^*$. The index $\hat{\pi}^* = 0.298$ suggests that the original table is far from independence, because about 30% of the population is estimated to be outside the model. A lower (approximate) 95% bound for $\pi^*$ is $\hat{\pi}_L^* = 0.236$.

Table 2 has $2nL^2 = 569.420$ on 12 degrees of freedom. The LRT statistic has extremely small $p$-values leading to rejection. The index $\pi^*$ is 0.104, and its lower limit is $\hat{\pi}_L^* = 0.091$. Given the potential for misclassification in either or both factors in Table 2, a misclassification rate of the order of 10% would not be surprising, and could explain the lack of fit. Regardless, the main conclusion is that the data in Table 2 are much closer to the hypothesis of row–column independence than are the data in Table 1.

5.2. *Kullback distance inference.* We next consider $K^2$ tubular model testing for our examples.

We solved the optimization problem (9) for a grid of Lagrange parameters $\pi$ chosen in $(0, 1)$. Note that our main purpose is to find the lower confidence limit of the tube radius. Thus we will need to compute for the tube radius as a function of $\pi$, that is $c = K^2(\pi \mathbf{d} + (1 - \pi)\mathbf{m}, \mathbf{m})$.

In Table 3, the LRT statistics $2nL^2$ along with $c$ are given for various values of $\pi$ for the data in Table 1. Note that when $c = 0$ the original model $\mathcal{M}$ is obtained. A monotone reduction in $L^2$-values is evident for increasing values of tubular radius $c$. For any value of $c \geq \hat{\rho}^*$, the tubular model will also be saturated yielding fit likelihood statistics of 0. A lower (approximate) 95% bound for $\rho^*$ is $\hat{\rho}^*_L = 0.101$, as noted in Table 3.

Table 4 gives analogous results for the data in Table 2. The value of $\hat{\rho}^*$ is about 0.011, and $\sqrt{\hat{\rho}^*} = 0.106$. The approximate 95% lower bound is $\hat{\rho}^*_L = 0.010$ with $\sqrt{\hat{\rho}^*_L} = 0.099$. Comparison of the latter with the value of 0.318 of Table 3 establishes the conclusion, just as for $\pi^*$, that the data in Table 2 are three times closer to the row–column independence than are the data in Table 1.

Note that the values of $\sqrt{\hat{\rho}^*}$, 0.369 and 0.106, are quite similar to those of $\hat{\pi}^*$, which were 0.298 and 0.104 in the two tables. This similarity between these (and other) measures gives some assistance for the overall interpretation of the distance, as one can select a distance with the most natural interpretation.

The close relationship between the index $\pi$ and the distances can effectively replace intensive grid search in $\pi$ for finding the values of $\pi$ that correspond to target values of $c$ and $2nL^2$. Indeed, both $\sqrt{K^2}$ and $\sqrt{L^2}$ behave nearly linearly in $\pi$ across along regions of $\pi$, reflecting a nearly Euclidean structure arising from the two distances together.

TABLE 3
*Likelihood statistics for the tubular model applied to Table 1*

| $\pi$ | $\sqrt{c}$ | $c$ | $2nL^2$ |
|---|---|---|---|
| 0.000 | 0.000 | 0.000 | 146.44 |
| 0.444 | 0.155 | 0.024 | 48.68 |
| 0.774 | 0.278 | 0.077 | 8.77 |
| **0.876** | **0.318** | **0.101** $(= \hat{\rho}^*_L)$ | **2.71** |
| 0.877 | 0.319 | 0.102 | 2.67 |
| 0.990 | 0.365 | 0.133 | 0.02 |
| **1.000** | **0.369** | **0.136** $(= \hat{\rho}^*)$ | **0.00** |

TABLE 4
*Likelihood statistics for the tubular model applied to Table 2*

| $\pi$ | $\sqrt{c}$ | $c$ | $2nL^2$ |
|---|---|---|---|
| 0.000 | 0.000 | 0.000 | 569.42 |
| 0.801 | 0.085 | 0.007 | 22.75 |
| 0.865 | 0.092 | 0.008 | 10.48 |
| **0.931** | **0.099** | **0.010** $(= \hat{\rho}^*_L)$ | **2.74** |
| 0.958 | 0.105 | 0.011 | 1.02 |
| **1.000** | **0.106** | **0.011** $(= \hat{\rho}^*)$ | **0.00** |

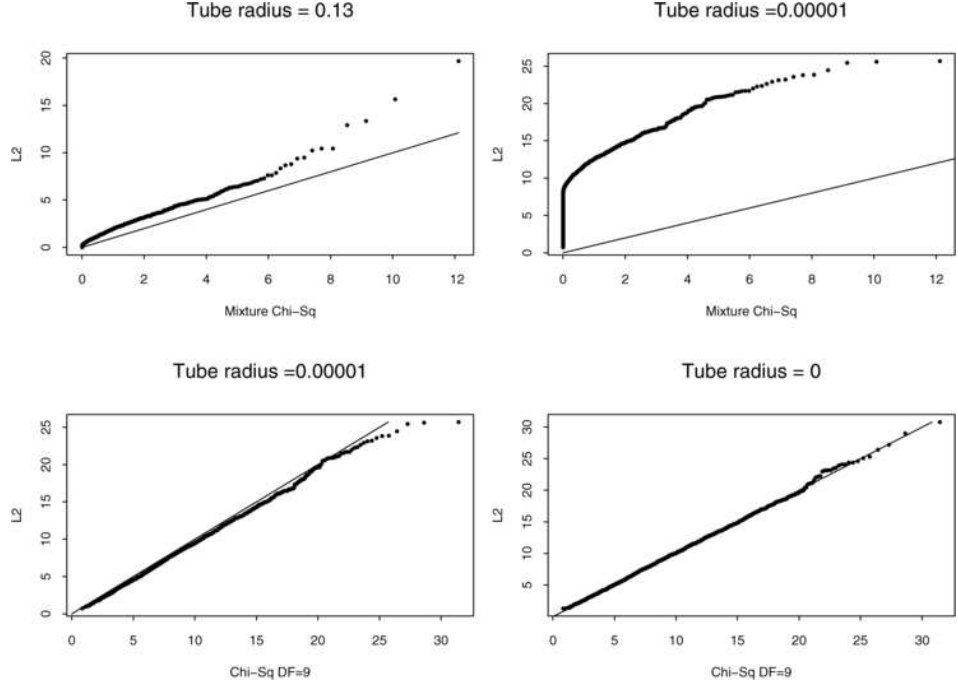

Fig. 1. *The simulated distribution of $L^2$ vs. theoretical asymptotic distribution.*

5.3. *Bootstrap results.* We here apply the bootstrap method to Table 1 ($4 \times 4$ table with $n = 592$ observations) and estimate the distribution of the LRT statistic $2nL^2$. Our bootstrap sample size is $B = 10^4$. We simulate $L^2$ for different tube radius, as $c$ takes different values from 0 to $\hat{\rho}^*$.

Asymptotic theory says that $2nL^2$ has a mixture distribution as $\frac{1}{2}\chi_0^2 + \frac{1}{2}\chi_1^2$ if $c > 0$, and $\chi_9^2$ if $c = 0$. Figure 1 gives the QQ plots that compare the simulated distributions with the asymptotic distributions of $L^2$ for different tubes. When $c = 0$ ($\pi = 0$), they match very well. But when $c$ is positive, the bootstrap and asymptotic distributions differs. Some of the discrepancy might be an effect of the curved boundary of the tube surface. When $c$ is positive but close to zero, the simulated and the asymptotic distributions strongly depart from each other, with the simulated sampling distribution much closer to $\chi_9^2$ than $\frac{1}{2}\chi_0^2 + \frac{1}{2}\chi_1^2$.

We calculated the upper 95% quantile of the bootstrap distribution to be used as the estimated critical value for testing $H_0$. Figure 2 plots the LRT statistics and the simulated 95% critical values for different tube radii. When $c > 0$, the simulated critical value and theoretical asymptotic critical value are different, especially for $c$ near zero.

The bootstrapped likelihood ratio test will reject tube models when the tube radius $c$ is smaller. In this example, we used linear interpolation and found when $\sqrt{c}$ is smaller than 0.300, the likelihood ratio test will reject $H_0$, so a lower 95% confidence limit for $\sqrt{K^2}$ is 0.300, compared to 0.318 using the asymptotic method. At this value of $\sqrt{c} = 0.300$ the simulated critical value for $2nL^2$ is 4.24, not very close to 2.70.

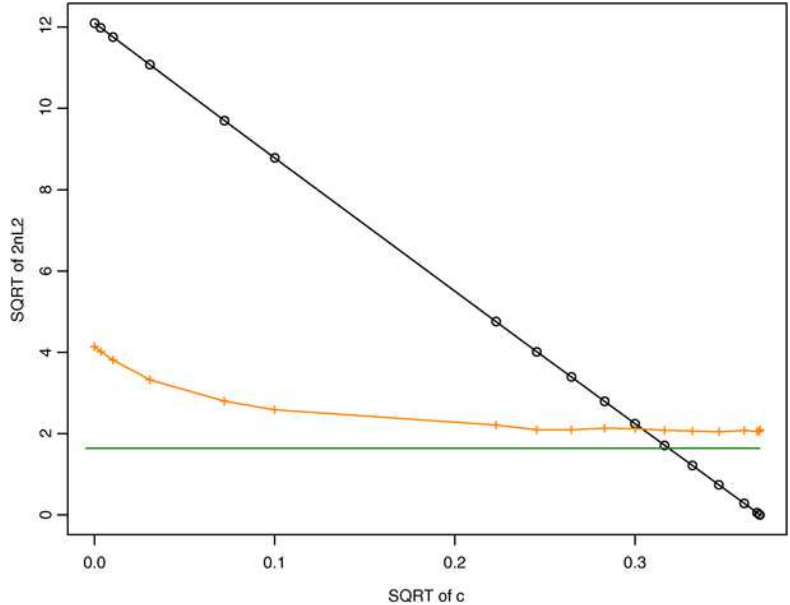

FIG. 2. *The LRT statistics ($\circ$) and the simulated critical values ($+$) against the tube radius.*

**6. Application to loglinear models.** In this section we will apply the tubular model to a multidimensional contingency table problem. We will compare our tubular index with the conventional LRT statistic for logistic regression models.

6.1. *The data and some results on loglinear models.* The data were collected from a sample survey of 8,036 army recruits, identified by race (C), geographic origin (R), location of current training camp (L) and their preference for camp location (P). The data of interest for this example is shown in Table 5 of Rudas (1991). Notice that the sample size is quite large, so one might find that some models do not fit, but still provide a high-quality approximation.

We will focus on comparison of various models. The candidate models are those with all main effects, or some two-order or three-order interactions.

The LRT and AIC, BIC indices, and our tubular model statistics for the candidate models are summarized in Table 6. Besides using the statistic

TABLE 5
*Preference of World War II Recruits for locations of training camp*

| Color (C) | Region (R) | Location (L) | Preference (P) | |
|---|---|---|---|---|
| | | | **North** | **South** |
| Black | North | North | 387 | 36 |
| | | South | 876 | 250 |
| | South | North | 383 | 270 |
| | | South | 381 | 1712 |
| White | North | North | 955 | 162 |
| | | South | 874 | 510 |
| | South | North | 104 | 176 |
| | | South | 91 | 869 |

$\hat{\rho}^*$, we calculate the lower confidence limit $\hat{\rho}_L^*$, by inverting the asymptotic critical value of the LRT statistic, with size $\alpha = 0.05$.

Except for the AIC and BIC values, which adjust for the effect of the number of parameters, the statistics $2nL^2$, $\hat{\rho}^*$ and $\hat{\rho}_L^*$ all improve as we add extra terms to models. For example, when Model 2 changes to Model 3, the goodness-of-fit improves. To minimize the AIC or BIC index, Model 4 would be selected as the best model among candidates. We have included AIC and BIC here for comparison, but note that they address the model selection problem from the point of view of risk, not distance. As discussed in Lindsay and Liu (2005), this makes it a measure of model quality that depends on the sample size $n$ used in the experiment.

6.2. *Interpretation of the tubular radius.* We then proceed to check the performance of the tubular statistics in assessing model fit. We note that the AIC and BIC select Model 4, which is also the most parsimonious model in which the lower confidence limit of the tubular index is zero. If we were to select a model by the rule "use the simplest model with the $\hat{\rho}_L^*$ value being zero", this would be equivalent to finding the smallest acceptable model by hypothesis testing. Using a model with $\hat{\rho}_L^* > 0$ would correspond to allowing a tolerance in the model fit. In order to make an acceptable lower limit for $\rho_L^*$, one needs a statistically meaningful interpretation.

TABLE 6
$2nL^2$, *AIC, BIC and tubular statistics for World War* II *Recruits data (k is the dimension of the parameter space)*

| Model | $k$ | $2nL^2$ | AIC | BIC | $\sqrt{\hat{\rho}^*}$ | $\sqrt{\hat{\rho}_L^*}$ |
|---|---|---|---|---|---|---|
| 1 Main effects | 4 | 4211.3 | 4219.3 | 4247.3 | 0.56 | 0.55 |
| 2 Model 1 + all 2-way interactions | 10 | 78.02 | 98.0 | 167.9 | 0.070 | 0.057 |
| 3 Model 2 + CRL | 11 | 24.96 | 46.9 | 123.9 | 0.039 | 0.026 |
| 4 Model 2 + CRL RLP | 12 | 1.45 | 25.5 | 109.4 | 0.0095 | 0 |
| 5 Model 2 + CRL CRP RLP | 13 | 0.68 | 26.7 | 117.6 | 0.0065 | 0 |
| 6 Model 2 + CRL CRP CLP RLP | 14 | 0.67 | 28.7 | 126.6 | 0.0064 | 0 |

In Table 7 we show the tubular index $\hat{\rho}^* = K^2(\mathbf{d}, \mathcal{M})$ for each model, and for companion, the mid-tube distance $T_{1/2}^2$, which is very similar in value after appropriate scaling, $(\rho^* \approx 4T_{1/2}^2)$. Although it is hard to give absolute meaning to $\sqrt{\rho^*}$, by comparing relative distance one can see improvement in fit if going from Model 1 to 2, but smaller gains thereafter.

Also shown in Table 7 is the value of the $\pi^*$ index. From this we note that for Model 3 only 6.5% of the population is not fit by the model. Relative magnitudes of $\sqrt{\hat{\rho}^*}$ are still similar to those of $\hat{\pi}^*$, but absolute magnitude is not as close as those in Section 5, again suggesting that Model 3 described a much larger fraction of the population.

In order to further interpret these distances, we use the index $N^*$ briefly described in Section 2. We implemented that idea here as follows. For each model $\mathcal{M}$, we used the size $\alpha = 0.05$ likelihood ratio test of $H_0 : \boldsymbol{\tau} \in \mathcal{M}$ against a general alternative. Given a nominal sample size $N$, we can then simulate samples from the empirical distribution $\mathbf{d}$. The proportion of rejections then provides us with an estimate of the power at sample size $N$ under the estimated alternative $\mathbf{d}$. We then find the sample size $N^*$ at which the power is 0.5. [For more details see Lindsay and Liu (2005)] This sample size then represents the maximum sample size at which samples from $\mathbf{d}$ are hard to distinguish from samples from $\mathcal{M}$.

In Table 7 the $N^*$ values are not reported for Models 4, 5 and 6 because these models were accepted in testing. That is, it is unreasonable to estimate the largest sample size for which the model is descriptive of the data when it is larger than the actual sample size. From this table we can see that the rejected Model 1, with root radius 0.56, is an extremely poor description of data samples of any size, larger than 5. Model 2, at root radius 0.070, describes samples of quite large size, 950. However, with small losses in parsimony, one could use Model 3 and gain considerable descriptive power over Model 2. Choice of a particular model in practice would then depend on the tradeoff between explanatory power and parsimony that one chooses to pursue.

TABLE 7
Tubular statistics, $\hat{\pi}^*$ and the credibility index $N^*$ for World War II recruits data

| Model | $k$ | $\sqrt{\hat{\rho}^*}$ | $\sqrt{4T_{1/2}^2}$ | $\hat{\pi}^*$ | $N^*$ |
|---|---|---|---|---|---|
| 1 | 4 | 0.563 | 0.517 | 0.715 | 5 |
| 2 | 10 | 0.0696 | 0.0696 | 0.143 | 950 |
| 3 | 11 | 0.0388 | 0.0394 | 0.065 | 3000 |
| 4 | 12 | 0.00945 | 0.00947 | 0 | $> n$ |
| 5 | 13 | 0.00650 | 0.00651 | 0 | $> n$ |
| 6 | 14 | 0.00642 | 0.00643 | 0 | $> n$ |

**7. Discussions.** In this paper we introduced the semiparametric tolerance region of a parametric model in the multinomial setting. We focused on the Kullback–Leibler distance to define the statistical tube. A likelihood ratio test procedure was developed for testing the tubular hypothesis that the true distribution is in such a tubular neighborhood of the model. The asymptotic and the simulated distributions of the LRT statistic were investigated and the lower confidence limit was then constructed.

We could use other distances to define the tube, such as Hellinger distance. The mid-tube distance proposed in Section 3.3 is another choice, since it is topologically equivalent to Hellinger distance. The most important feature relative to tubal inference is that when one uses Hellinger distance both for the tube distance and for the test statistic, one gets a closed form solution for the single element model as in Theorem 1 [see Liu (2003)]. One can mimic much of the development of this paper including fast algorithms.

Finally, the concept of tubular tolerance region could be extended to continuous distributions using the empirical likelihood. So far we applied our methods only for multinomial models. The LRT approach in the models can be replaced by empirical likelihood methods, enabling one to generalize the concept of tube to continuous distributions. This is under investigation by the authors.

Department of Mathematics and Statistics  
Georgia State University  
Atlanta, Georgia 30022  
USA  
E-mail: matjxl@langate.gsu.edu

Department of Statistics  
Pennsylvania State University  
University Park, Pennsylvania 16802  
USA  
E-mail: bgl@psu.edu